\documentclass[a4paper]{amsart}

\usepackage{amscd,amsthm,amsfonts,amssymb,amsmath}
\usepackage{mathrsfs}
\usepackage[colorlinks=true]{hyperref}
\usepackage{fullpage}
\usepackage[all]{xy}
\usepackage{xcolor}

    \newcommand{\BA}{{\mathbb {A}}} 
    \newcommand{\BC}{{\mathbb {C}}} 
    \newcommand{\BE}{{\mathbb {E}}} \newcommand{\BF}{{\mathbb {F}}}

     \newcommand{\BL}{{\mathbb {L}}}

    \newcommand{\BQ}{{\mathbb {Q}}}

     \newcommand{\BZ}{{\mathbb {Z}}}

     \newcommand{\CB}{{\mathcal {B}}}

    \newcommand{\CO}{{\mathcal {O}}}

     \newcommand{\RN}{{\mathrm {N}}}

     \newcommand{\fp}{{\mathfrak{p}}}

     \newcommand{\fA}{{\mathfrak{A}}}

    \newcommand{\Char}{{\mathrm{Char}}}

    \newcommand{\Gal}{{\mathrm{Gal}}} \newcommand{\GL}{{\mathrm{GL}}}
    
    \newcommand{\Hom}{{\mathrm{Hom}}}

    \newcommand{\Ind}{{\mathrm{Ind}}}

    \renewcommand{\mod}{\ \mathrm{mod}\ }

      \newcommand{\Supp}{{\mathrm{Supp}}}

    \newcommand{\Vol}{{\mathrm{Vol}}}\newcommand{\vol}{{\mathrm{vol}}}

        \newcommand{\Tr}{\mathrm{Tr}}

\newcommand{\matrixx}[4]{\begin{pmatrix}
#1 & #2 \\ #3 & #4
\end{pmatrix} }        

    \newcommand{\pair}[1]{\langle {#1} \rangle}

    \newcommand{\ov}{\overline}

    \newcommand{\ra}{\rightarrow}

\newcommand{\lrb}[1]{\left(#1\right)}

\newcommand{\B}{\ensuremath{{\mathbb{B}}}}
\newcommand{\C}{\ensuremath{{\mathbb{C}}}}

\newcommand{\F}{\ensuremath{{\mathbb{F}}}}

 \newcommand{\BetaI}[1]{{\left\{#1 \right\} }}
\newcommand{\E}{\ensuremath{{\mathbb{E}}}}
\renewcommand{\L}{\ensuremath{{\mathbb{L}}}}
\newcommand{\Nm}{\ensuremath{\text{Nm}}}
\newcommand{\WaldsI}{I(\varphi,\chi)}
\newcommand{\zxz}[4]{\begin{pmatrix} #1 & #2 \\ #3 & #4 \end{pmatrix}}

\newcommand{\Cor}[1]{}
\newcommand{\New}[1]{#1}

    \theoremstyle{plain}
    \newtheorem{thm}{Theorem}[section] \newtheorem{coro}[thm]{Corollary}
    \newtheorem{lem}[thm]{Lemma}  \newtheorem{prop}[thm]{Proposition}
     \newtheorem{defn}[thm]{Definition}

\theoremstyle{remark} \newtheorem{remark}{Remark}[section]
\theoremstyle{remark} 
\theoremstyle{remark} 

    \numberwithin{equation}{section}

\begin{document}
\title[Waldspurger's period integral for newforms]{Waldspurger's period integral for newforms}
\author[Yueke Hu, Jie Shu and Hongbo Yin]{Yueke Hu, Jie Shu and Hongbo Yin}

\begin{abstract}
In this paper we discuss Waldspurger's local period integral for newforms in new cases. The main ingredient is the work \cite{HN18} on Waldspurger's period integral using the minimal vectors, and the explicit relation between the newforms and the minimal vectors. We use a representation theoretical trick to simplify computations for newforms. As an example, we compute the local integral coming from a special arithmetic setting which was used to study 3-part full BSD conjecture in \cite{HSY}.
\end{abstract}
\address{Department of Mathematics, ETH, Zurich, Switzerland}
\email{huyueke2012@gmail.com}

\address{School of Mathematical Sciences, Tongji University, Shanghai 200092}
\email{shujie@tongji.edu.cn}
\address{School of Mathematics, Shandong University, Jinan 250100,
P.R.China}
\email{yhb2004@mail.sdu.edu.cn}


\thanks{Yueke Hu is supported by SNF-169247; Jie Shu is supported by NSFC-11701092; Hongbo Yin is supported by NSFC-11701548 and The Fundamental Research Funds of Shandong University.}

\maketitle


\section{Introduction}
In this paper we study Waldspurger's local period integral for newforms in new cases. Most of the paper is purely local, so we shall skip subscript $v$ when there is no confusion.

Suppose that $p\neq 2$. Let $\pi$ be a smooth irreducible self-dual representation of $\GL_2$ over a $p$-adic field $\BF$ with central character $w_\pi$, and $\chi$ be any character over a quadratic field extension $\BE/\BF$ such that  $w_\pi \chi|_{\BF^\times}=1$. 
For any test vector $\varphi\in\pi$, denote the local Waldspurger's period integral against a character $\chi$ on $\BE^\times$ by
\begin{equation}
\WaldsI=\{\varphi,\varphi\}\text{\ where } \{\varphi,\varphi'\}=\int\limits_{t\in \BF^\times\backslash \BE^\times} \Phi_{\varphi,\varphi'}(t)\chi(t)dt.
\end{equation}
Here $\Phi_{\varphi,\varphi'}(t)=(\pi(t)\varphi,\varphi')$ is the matrix coefficient associated to $\varphi$, $\varphi'$. 
We are particularly concerned about the case when $\pi$ is a supercuspidal representation and $\varphi$ is a proper translate of the newform for a fixed embedding of $\BE$ (or equivalently $\varphi$ is a fixed newform and the embedding of $\BE$ can vary).

The explicit knowledges of the local period integral and the particular choice of test vectors have been useful for analytic number theory as well as arithmetic geometry problems. For example, they were used to study the moments and the subconvexity bound of L-functions in \cite{FileMartinPitale:17} and \cite{Wu:16a}.
In \cite{CST14}, they were used to give general explicit Gross-Zagier and Waldspurger formulae, which is important to attack the refined BSD conjecture. 
In the most recent example, we used the explicit Gross-Zagier formula in \cite{HSY} to study the 3-part full BSD conjecture for the elliptic curves related to the Sylvester conjecture. The byproduct of this paper  is a special example of  the local period integral which is not covered in the previous literatures, and is used in \cite{HSY} to establish the explicit Gross-Zagier formula there.
\subsection{A short history of local test vector}
The study of test vectors when there are ramifications was initiated in \cite{GrossPrasad:91a}. It assumes disjoint ramifications, and describes a test vector in terms of invariance by proper compact subgroups. 
The work in \cite{FileMartinPitale:17}   gives test vectors in more general situations on $\GL_2$ side. In particular it solves the case when $\E$ is split. When $\E$ is a quadratic field extension, 
its method can deal with the range $c(\pi_{\chi})>c(\pi)$, (here $\pi_{\chi}$ is the representation of $\GL_2$ associated to $\chi$ via the Theta correspondence or Langlands correspondence,) and the test vectors used are diagonal translates of newforms.

The recent work \cite{HN18} under mild assumptions provided test vectors for the complementary range $c(\pi_{\chi})\leq c(\pi)$ when $\pi$ is supercuspidal, using a new type of test vectors called minimal vectors. A particular minimal vector is given in the Kirillov model in Lemma \ref{lem:toricnewforminKirillov} below, and any single translate of this element is still considered to be a minimal vector.  Such test vectors arise naturally from compact induction theory for supercuspidal representations and have better properties than the standard newforms in terms of their matrix coefficients and Whittaker functionals.   For some applications however (in particular with \cite{HSY} in mind), it's still necessary to understand the local period integral for the classical newforms. The purpose of this paper is thus to make use of \cite{HN18} to predict the local period integral for the newforms.
\subsection{Main results}
The goal of this paper is mainly to present the method, rather than exhausting all the cases. So we restrict ourselves to a special setting to avoid lengthy discussions. We shall assume that $\pi$ is associated to a character $\theta$ over a ramified extension $\BL$ via compact induction theory with $c(\theta)=2n$, $\theta|_{\F^\times}=w_\pi=1$, 
 $\BE\simeq\BL$ is also ramified,   and $c(\pi_{\chi})\leq c(\pi)$. We also assume $\epsilon(\pi_\BE\times\chi)=1$ so that $\Hom_{\BE^\times}(\pi,\chi^{-1})\neq 0$.

For simplicity we shall fix the embedding of $\BL$ into $\mathrm{M}_2$ by
\begin{equation}
   x+y\sqrt{D'}\mapsto \matrixx{x}{y}{yD'}{x},
\end{equation}
where $v(D')=1$. Any other embedding differs by a conjugation, which is effectively equivalent to a single translate for the test vector.

Using the explicit Kirillov model for a particular minimal vector $\varphi_0$ in Lemma \ref{lem:toricnewforminKirillov}, we can write the newform as a sum of minimal vectors in Corollary \ref{Cor:RelationNewMinimal}. Correspondingly we can write the period integral for the newform as a sum of period integrals for the minimal vectors in Corollary \ref{Cor:PeriodIntRelation}, which specialises in our setting to the following
\begin{equation}\label{Eq:Intro-Periodrelation}
\{\widetilde{\varphi_{new}},\widetilde{\varphi_{new}}\}=\frac{1}{(q-1)q^{\lceil \frac{n}{2 }\rceil-1}}\sum\limits_{x,x'\in (\CO_\BF/\varpi^{\lceil \frac{n}{2}\rceil}\CO_\BF)^\times}\{\varphi_x,\varphi_{x'}\}.
\end{equation}
Here $\widetilde{\varphi_{new}}=\pi\lrb{\zxz{\varpi^n}{0}{0}{1}}\varphi_{new}$ is a diagonal  translate of the standard newform, and $\varphi_x=\pi\lrb{\zxz{x}{0}{0}{1}}\varphi_0$. The work in \cite{HN18} computed the diagonal terms $\{\varphi_x,\varphi_{x}\}$ which is either 0 or some constant only depending on the associated conductors.

In the case when there is a single diagonal term $\{\varphi_x,\varphi_{x}\}$ in the family which is nonvanishing, one can use the representation theoretical trick in Lemma \ref{lem:Gl2-newform-crossterms} to easily get the following
\begin{prop}[Proposition \ref{Prop:singlenonvanishingcase}]
When $\BetaI{\varphi_x,\varphi_x}\neq 0$ for a single $x\in (\CO_\BF/\varpi^{\lceil \frac{n}{2}\rceil}\CO_\BF)^\times$, all off-diagonal terms vanish and
\begin{equation}
 \BetaI{\widetilde{\varphi_{new}},\widetilde{\varphi_{new}}}=\frac{1}{(q-1)q^{\lceil \frac{n}{2}\rceil-1}}\BetaI{\varphi_x,\varphi_{x}}.
\end{equation}
\end{prop}
Note here that $\BetaI{\varphi_x,\varphi_{x}}$ is already computed in \cite{HN18}.

The more challenging case is when there are several  non-vanishing diagonal terms. The corresponding off-diagonal terms will also be non-vanishing and have the same absolute value as the diagonal terms by Lemma \ref{lem:Gl2-newform-crossterms}. The main innovation of this paper is that we devised a way to detect the phase factor for the off-diagonal terms with relatively simple computations. 
Roughly speaking, the sizes of the non-vanishing off-diagonal terms are the same as those for the non-vanishing diagonal terms. The support of the integral (see Definition \ref{Defn:suppofint}) can also be computed directly. By comparing the support of the integral with the size of the integral, we shall see that the integrand must be constant on the support of the integral. This constant is the phase factor and can be easily detected by taking special values.

In particular we obtain the following result.
\begin{prop}[Proposition \ref{Prop:TestingonNewform}]\label{Prop:Intromain}
Suppose that $\BE\simeq\BL$ are ramified, $c(\theta)=2n$, $w_\pi=1$, $\epsilon(\pi_\BE\times\chi)=1$, and $0<c(\theta\overline{\chi})=2l\leq 2n$ with $n-l$ even. Then
$$I(\widetilde{\varphi_{new}},\chi)=\frac{1}{(q-1)q^{\lceil \frac{n}{2 }\rceil-1}}\frac{1}{q^{\lfloor l/2\rfloor}} (1+\theta\chi(\sqrt{D}))^2    .$$
In particular it is either $0$ or asymptotically $\frac{1}{C(\pi\times\pi_{\chi})^{1/4}}$, where $C(\pi\times\pi_{\chi})$ is the conductor for the associated Rankin-Selberg $L-$function.
\end{prop}
\begin{remark}
\begin{enumerate}
\item One feature used in our particular setting in Proposition \ref{Prop:Intromain} is that there are exactly two non-vanishing diagonal terms in the expansion \eqref{Eq:Intro-Periodrelation}. According to \cite{HN18}, this is true in most of the other cases.
\item Unlike the range $c(\pi_\chi)>c(\pi)$, we see here the additional obstruction for using diagonal translate of newforms. 
\item Hopefully the same strategy works for more general cases, at least for those mentioned in (1).
\end{enumerate}

\end{remark}

We would like to point out that both Proposition \ref{Prop:singlenonvanishingcase} and Proposition \ref{Prop:TestingonNewform} are new cases not known before.
As a motivation for these discussions, we give the local period integral in Proposition \ref{Prop:TestingForNew} coming from a special arithmetic setting, 
where Proposition \ref{Prop:singlenonvanishingcase} alone turns out to be enough. Such result is one of the key ingredients used in \cite{HSY} to discuss the 3-part full BSD conjecture. 

\section{Notations and preliminary results}\label{minimal vector}
\subsection{Notations and basics}\label{prelimi}



For a real number $a$, let $\lfloor a\rfloor\leq a$ be the largest possible integer, and $\lceil a\rceil \geq a$ be the smallest possible integer.

Let $\F$ be a $p$-adic field with \New{residue field of order} $q$, uniformizer $\varpi=\varpi_\F$, ring of integers $\CO_\F$ and $p$-adic valuation $v_\F$. Let $\psi$ be an additive character of $\F$. 
Assume that $2\nmid q$. For $n\geq 1$, let $U_\F(n)=1+\varpi_\F^n \CO_\F$.
Let $\pi$ be a supercuspidal representation over $\F$ with central character $w_\pi=1$. 

Let $\L$  be a quadratic field extension over $\F$. Let $e_\L=e(\L/\F)$ be the ramification index and $v_\L$ be the valuation on $\L$. Let $\varpi_\L$ be a uniformizer for $\L$. When $\L$ is unramified we shall identify $\varpi_\L$ with $\varpi_\F$. Otherwise we suppose that $\varpi_\L^2=\varpi_\F$. Let $x\mapsto \overline{x}$ to be the unique nontrivial involution of $\L/\F$.
 Let $\psi_\L=\psi\circ \Tr_{\L/\F}$. One can make similar definitions for a possibly different quadratic field extension $\E$. Note that we shall assume that $\varpi_\E^2=\xi\varpi_\F$ for $\xi\in \CO_\F^\times-(\CO_\BF^\times)^2$ if $\E$ and $\L$ are both ramified and distinct.

For $\chi$ a multiplicative character on $\F^\times$, let $c(\chi)$ be  the smallest integer such that $\chi$ is trivial on $1+\varpi_\F^{c(\chi)} \CO_\F$. Similarly $c(\psi)$ is the smallest integer such that  $\psi$ is trivial on $\varpi_\F^{c(\psi)}\CO_\F$. We choose $\psi$  to be unramified, or equivalently,   $c(\psi)=0$. Then $c(\psi_\L)=-e_\L+1$.
Let $c(\pi)$ be the power of the conductor of $\pi$.

When $\chi$ is a character over a quadratic extension, denote $\overline{\chi}(x)=\chi(\overline{x})$.
\Cor{\begin{lem}
Suppose that $p$ is large enough.
For a multiplicative character $\nu$ over $\F$ with $c(\nu)\geq 2$, there exists $\alpha_\nu\in \F^\times$ with $v_\F(\alpha_\nu)=-c(\nu)+c(\psi)$ such that for $u\in \varpi_\F \CO_\F$
\begin{equation}
\nu(1+u)=\psi(\alpha_\nu \log(1+u)) 
\end{equation}
where $\log(1+u)$ is the standard Taylor expansion for logarithm 
\begin{equation}
\log(1+u)=u-\frac{u^2}{2}+\frac{u^3}{3}+\cdots.
\end{equation}
In particular when $p>2$ we have that 
\begin{equation}
\nu(1+u)=\psi(\alpha_\nu (u-\frac{u^2}{2})) \end{equation}
for any $u\in \varpi_\F^{\lfloor c(\nu)/2\rfloor} \CO_\F$.
\end{lem}}
\New{\begin{lem}\label{Lem:DualLiealgForChar}
For a multiplicative character $\nu$ over $\F$ with $c(\nu)\geq 2$, there exists $\alpha_\nu\in \F^\times$ with $v_\F(\alpha_\nu)=-c(\nu)+c(\psi)$ such that
\begin{equation}\label{eq:alphatheta}
\nu(1+u)=\psi(\alpha_\nu u) \end{equation}
for any $u\in \varpi_\F^{\lceil c(\nu)/2\rceil} \CO_\F$. $\alpha_\nu$ is determined $\mod \varpi_\F^{-\lceil c(\nu)/2\rceil+c(\psi)} \CO_\F$.
\end{lem}
One can easily check this lemma by using that $\nu(1+u)$ becomes an additive character in $u$ for $u\in \varpi_\F^{\lceil c(\nu)/2\rceil} \CO_\F$.
}

We shall also need some basic results for compact induction theory and matrix coefficient.
In general, let $G$ be a unimodular locally profinite group with center $Z$. Let $H\subset G$ be an open and closed subgroup containing $Z$ with $H/Z$ compact. Let $\rho$ be an irreducible smooth representation of $H$ with unitary central character and 
\New{$$\pi=c-\Ind_H^G(\rho)=\{f:G\rightarrow \rho\mid f(hg)=\rho(h)f(g)\ \forall\ h\in H, \text{\ $f$ is compactly supported}\}.$$}
  By the assumption on $H/Z$, $\rho$ is automatically unitarisable, and we shall denote the unitary pairing on $\rho$ by $\pair{\cdot,\cdot}_{\rho}$. Then one can define the unitary pairing on $\pi$ by
\begin{equation}
\pair{\phi,\psi}=\sum\limits_{x\in H\backslash G}\pair{\phi(x),\psi(x)}_{\rho}.
\end{equation}
If we let $y\in H\backslash G$ and $\{v_i\}$ be a basis for $\rho$, the elements
\begin{equation}\label{Eq:CompactIndBasis}
 f_{y,v_i}(g)=\begin{cases}
\rho(h)v_i,&\text{\ if \ }g=hy\in Hy;\\
0,&\text{\ otherwise.}
\end{cases}
\end{equation} 
form a basis for $\pi$.
\begin{lem}\label{lem3.3:matrixcoeffInduction}
For $y,z\in H\backslash G$,
\begin{equation}
\pair{\pi(g)f_{y,v_i},f_{z,v_j}}=\begin{cases}
\pair{\rho(h)v_i,v_j}_{\rho}, &\text{\ if\ }g=z^{-1}hy\in z^{-1}Hy;\\
0,&\text{\ otherwise}.
\end{cases}
\end{equation}
\end{lem}

We also recall the multiplicity one result and Tunnell-Saito's epsilon value test.
\begin{thm}[\cite{Tunnell:83a}\cite{Saito:93a} ]\label{Tunnell}
Suppose that $w_\pi=\chi|_{\F^\times }$. Let $\pi^\B$ be the image of $\pi$ under the Jacquet-Langlands correspondence.
The space $\Hom_{\E^\times }(\pi^\B\otimes \chi^{-1},\C)$ is at most one-dimensional. It is nonzero if and only if 
\begin{equation}
\epsilon(\pi_\E\times\chi^{-1})=\chi(-1)\epsilon(\B).
\end{equation}
Here $\pi_\E$ is the base change of $\pi$ to $\E$. $\epsilon(\B)=1$ if it is a matrix algebra, and $-1$ if it's a division algebra.
\end{thm}

\subsection{Compact induction theory for supercuspidal representations and minimal vectors}
We shall review the compact induction theory for supercuspidal representations on $\GL_2$. For more details, see \cite{BushnellHenniart:06a}.

We shall fix the embeddings and work out everything explicitly. For $v_\F(D')=0$ or $1$, we shall refer to the following embedding of a quadratic field extension $\BL=\BF(\sqrt{D'})$ as a standard embedding:
\begin{equation}\label{Eq:standardembeddingL}
   x+y\sqrt{D'}\mapsto \matrixx{x}{y}{yD'}{x}.
\end{equation}

The supercuspidal representations are parametrized via compact induction by characters $\theta$ over some quadratic field extension $\L$. More specifically we have the following quick guide. \begin{enumerate}
\item[Case 1.]$c(\pi)=2n+1$ corresponds to $e_\BL=2$ and $c(\theta)=2n$ .
\item[Case 2.] $c(\pi)=4n$ corresponds to $e_\BL=1$ and $c(\theta)=2n$.
\item[Case 3.] $c(\pi)=4n+2$ corresponds to $e_\BL=1$ and $c(\theta)=2n+1$ .
\end{enumerate}

\begin{defn}
For $e_{\BL}=1,2$, let $$\fA_{e_{\BL}}=\begin{cases}
M_2(\CO_{\BF}), \text{\ if }e_\BL=1,\\
\matrixx{ \CO_\BF}{\CO_\BF}{\varpi \CO_\BF}{ \CO_\BF},\text{\ otherwise}.
\end{cases}
$$
Its Jacobson radical is given by
$$\CB_{e_\BL}=\begin{cases}
\varpi M_2(\CO_{\BF}), \text{\ if }e_\BL=1,\\
\matrixx{\varpi \CO_\BF}{\CO_\BF}{\varpi \CO_\BF}{\varpi \CO_\BF},\text{\ otherwise}.
\end{cases}$$
\end{defn}

Define the filtrations of compact open subgroups as follows:
\begin{equation}
K_{\fA_{e_\BL}}(n)=1+\CB_{e_\BL}^n,\ \ \  U_{\BL}(n)=1+\varpi_{\BL}^n\CO_{\BL}.
\end{equation}
Note that each $K_{\fA_{e_\BL}}(n)$ is normalised by $\BL^\times$ which is embedded as in \eqref{Eq:standardembeddingL}.

Denote 
$J=\BL^\times K_{\fA_{e_\BL}}(\lfloor c(\theta)/2\rfloor)$, $J^1=U_\BL(1)K_{\fA_{e_\BL}}(\lfloor c(\theta)/2\rfloor)$, $H^1=U_\BL(1)K_{\fA_{e_\BL}}(\lceil c(\theta)/2\rceil)$. Then $\theta$ on $\BL^\times$ can be extended to be a character $\tilde{\theta}$ on $H^1$ by
\begin{equation}\label{Eq:thetatilde}
\tilde{\theta}(l (1+x))=\theta(l)\psi\circ\Tr (\alpha_\theta x),
\end{equation}
where $l\in \BL^\times$, $1+x\in K_{\fA_{e_\BL}}(\lceil c(\theta)/2\rceil)$ and $\alpha_\theta\in \BL^\times\subset M_2(\BF)$ is associated to $\theta$ by Lemma \ref{Lem:DualLiealgForChar} under the fixed embedding.

When $c(\theta)$ is even, $H^1=J^1$ and $\tilde{\theta}$ can be further extended to $J$ by the same formula. In this case denote $\Lambda=\tilde{\theta}$  and $\pi=c-\Ind_J^G\Lambda$.

When $c(\theta)$ is odd, $J^1/H^1$ is a two dimensional vector space over the residue field. This case only occurs when $c(\pi)=4n+2$ as listed above.
Then there exists a $q-$dimensional representation $\Lambda$ of $J$ such that
$\Lambda|_{H^1}$ is a multiple of $\tilde{\theta}$, and
$\Lambda|_{\BL^\times}=\oplus \theta\nu$ where $\nu$ is over $\L^\times$, $c(\nu)=1$ and $\nu|_{\BF^\times}=1$. More specifically, let $B^1$ be any intermediate group between $J^1$ and $H^1$ such that $B^1/H^1$ gives a polarisation of $J^1/H^1$ under the  pairing given by \Cor{$\psi_{\alpha_\theta}$.}
\New{$$(1+x,1+y)\mapsto \psi\circ \Tr (\alpha_\theta [x, y]).$$}
 Then $\tilde{\theta}$ can be extended to $B^1$ by the same formula \eqref{Eq:thetatilde} and  $\Lambda|_{J^1}=\Ind_{B^1}^{J^1}\tilde{\theta}$.
 Again $\pi=c-\Ind_{J}^G\Lambda$ in this case.

 In the case $J^1=H^1$, we take $B^1=J$ for uniformity. In either cases, we have $w_\pi=\theta|_{\BF^\times}$.

\begin{defn}
There exists a unique element $\varphi_0\in \pi$ such that $B^1$ acts on it by $\tilde{\theta}$. (Type 1 minimal vector in the notation of \cite{HN18}.)
We also call any single translate $\pi(g)\varphi_0$ a minimal vector,
\end{defn}
Note that the conjugated group $gB^1g^{-1}$ acts on $\pi(g)\varphi_0$ by the conjugated character $\tilde{\theta}^g$ and They form an orthonormal basis of $\pi$.
\begin{coro}\label{Cor:MCofGeneralMinimalVec}
Let $\Phi_{\varphi_0}$ be the matrix coefficient associated to a minimal vector $\varphi_0$ as above. Then $\Phi_{\varphi_0}$ is supported on $J$, and
\begin{equation}
\Phi_{\varphi_0}(bx)=\Phi_{\varphi_0}(xb)=\tilde{\theta}(b)\Phi_{\varphi_0}(x)
\end{equation}
for any $b\in B^1$. Furthermore when $\dim \Lambda\neq 1$, $\Phi_{\varphi_0}|_{J^1}$ is supported only on $B^1$.
\end{coro}
Note that $\varphi_0$ is basically $f_{1,v_i}$ as in \eqref{Eq:CompactIndBasis} for the coset representative $1\in J\backslash G$. The corollary \New{follows immediately} from the definition of $\varphi_0$ and Lemma \ref{lem3.3:matrixcoeffInduction}.

\subsection{Local Langlands correspondence and compact induction}\label{SubSec:Langlands-compactInd}
Here we describe the relation between the compact induction parametrisation and the local Langlands correspondence. See \cite{BushnellHenniart:06a} Section 34 for more details.

For a field extension $\L/\F$ and an additive character $\psi$ over $\F$, let $\lambda_{\L/\F}(\psi)$ be the Langlands $\lambda-$function in \cite{Langlands}. When $\L/\F$ is a quadratic field extension, let $\eta_{\L/\F}$ be the associated quadratic character. By \cite[Lemma 5.1]{Langlands}, we have for $\psi_\beta(x)=\psi(\beta x)$,
\begin{equation}\label{Eq:LanglandsLambdaFun}
\lambda_{\L/\F}(\psi_\beta)=\eta_{\L/\F}(\beta)\lambda_{\L/\F}(\psi).
\end{equation}
\begin{defn}\label{Defn:DeltaTheta}

\begin{enumerate}
\item If $\BL$ is inert, define $\Delta_\theta$ to be the unique unramified character of $\BL^\times$ of order 2.

\item If $\BL$ is ramified and $\theta$ is a character over $\BL$ with $c(\theta)>0$ even, associate $\alpha_\theta$ to $\theta$ as in Lemma \ref{Lem:DualLiealgForChar}.
Then define $\Delta_\theta$ to be the unique level 1 character of $\BL^\times$ such that
\begin{align}
\Delta_\theta|_{\BF^\times}=\eta_{\BL/\BF}, \Delta_\theta(\varpi_\BL)=\eta_{\BL/\BF}(\varpi_\BL^{c(\theta)-1}\alpha_\theta)\lambda_{\BL/\BF}^{c(\theta)-1}(\psi).
\end{align}
\end{enumerate}

\end{defn}
Note that in  \cite{BushnellHenniart:06a}   $\psi$ is chosen to be of level 1. We have adapted the formula there to our choice of $\psi$ using \eqref{Eq:LanglandsLambdaFun}.
The definition is also independent of the choice of $\varpi_\BL$.    
\begin{thm}\label{LLC}
If $\pi$ is associated by compact induction to a character $\theta$ over a quadratic extension $\BL$, then its associated Deligne-Weil representation by local Langlands correspondence is $\sigma=\Ind_{\BL}^{\BF} (\Theta)$, where $\Theta=\theta\Delta_\theta^{-1}$, or equivalently $\theta=\Theta\Delta_{\Theta}$.
\end{thm}
Note here that $\Theta$ and $\theta$ always differ by a  level $\leq 1$ character, so $\alpha_\Theta$ \New{can be chosen to be the same as} $\alpha_\theta$ in Lemma \ref{Lem:DualLiealgForChar} and $\Delta_\Theta=\Delta_\theta$.

\subsection{Using minimal vectors for Waldspurger's period integral }\label{Sec:WaldsForMinimalVec}
\New{Now we review the local Waldspurger's period integral for minimal vectors, in the setting $e_\L=2$ and $\E\simeq \L$. For more details and proofs, see \cite{HN18}, and in particular its appendix for explicit treatment.}

For simplicity we pick $$D'=\frac{1}{\alpha_\theta^2\varpi_\L^{2c(\theta)}},$$
identify $\frac{1}{\alpha_\theta\varpi_\L^{c(\theta)}}$ with $\sqrt{D'}$ and $\L$ with $\F(\sqrt{D'})$, and use the standard embedding \eqref{Eq:standardembeddingL}.

By this choice, we have $v_\F(D')=0$ if $e_\L=1$ and $v_\F(D')=1$ if $e_\L=2$.
\begin{equation}\label{eq2.1:specialembedding}
\alpha_\theta= \frac{1}{\varpi_\L^{c(\theta)}}\frac{1}{\sqrt{D'}} \mapsto \frac{1}{\varpi^{c(\theta)/e_\L}}\zxz{0}{\frac{1}{D'}}{1}{0}.
\end{equation}
Such choice is not essential. A different choice will result in slightly different formulation, for example, in \eqref{eq:new-necessary-same-ramified}, but the final results are similar.

Note that when $e_\L=2$, $c(\theta)$ must be even when $\theta|_{\F^*}=1$. 
 Assume that $\BE=\BF(\sqrt{D})$ for $v_\F(D)=0,1$ is also embedded as
\begin{equation}\label{Eq:standardembedding}
   x+y\sqrt{D}\mapsto \matrixx{x}{y}{yD}{x}. 
\end{equation}

In \cite{HN18}, test vectors of form $\pi(g)\varphi_0$ were used to study $\WaldsI$ for general combinations of $\pi$, $\BE$ and $\chi$. 
For the purpose of this paper we shall only review the case when $\BL\simeq\BE$ are ramified,  $w_\pi=1$, $c(\pi)=2n+1$ is odd and $c(\pi_{\chi})\leq c(\pi)$. 
We shall normalize the Haar measure on $\BE^\times$ so that $\Vol( \CO_\BF^\times\backslash \CO_\BE^\times)=1$. Then $\Vol(\BF^\times\backslash \BE^\times)=2$.

In this case we use test vectors of the form  $$\varphi=\pi\lrb{\matrixx{1}{u}{0}{1}\matrixx{v}{0}{0}{1}}\varphi_0$$
for some $u\in \CO_\BF,v\in \CO_\BF^\times$. Recall that when $e_\BL=2$, 
$$K_{\fA_{e_\BL}}(n)=1+\matrixx{\varpi^{ \lceil n/2\rceil}\CO_\BF}{\varpi^{ \lfloor n/2\rfloor}\CO_\BF}{\varpi^{ \lfloor n/2\rfloor+1}\CO_\BF}{\varpi^{ \lceil n/2\rceil}\CO_\BF},$$
and $J=\BL^\times K_{\fA_{e_\BL}}(n)$ acts on $\varphi_0$ by a character, so we can assume that $v\in (\CO_\BF/ \varpi^{\lceil n/2\rceil}\CO_\BF)^\times$.

There are two situations depending on $\min \{c(\theta\chi), c(\theta\overline{\chi})
\}$. Note that for the embedding of $\BE$ fixed above, we have
\begin{equation}
\zxz{-1}{0}{0}{1} \matrixx{x}{y}{yD}{x}\zxz{-1}{0}{0}{1}=\matrixx{x}{-y}{-yD}{x}
\end{equation}
and thus
\begin{equation}
I(\varphi,\chi)=I\lrb{\pi\lrb{\zxz{-1}{0}{0}{1}}\varphi, \overline{\chi}}.
\end{equation}
So we shall always assume that $c(\theta{\overline{\chi}})\leq c(\theta\chi)$.

The first situation is when $c(\theta{\overline{\chi}})=0$, then the Tunnell-Saito's test requires $\theta{\overline{\chi}}$ to be trivial for $\WaldsI$ to be nonzero. In that case we can take $u=0$ and then there is a unique $v\mod \varpi^{\lceil n/2\rceil}$ such that $\WaldsI\neq 0$, and 
for this $v$ we have \begin{equation}
 \WaldsI\New{=\vol(\BF^\times\backslash\BE^\times)}=2.
\end{equation}

The second situation is when $0<c(\theta{\overline{\chi}})=2l\leq 2n$. In this case, $\alpha_{\theta{\overline{\chi}}}$ can be associated to $\theta\overline{\chi}$ by Lemma \ref{Lem:DualLiealgForChar} and $\varphi$ would be a test vector if $v$, $u$ are solutions of the following quadratic equation:
\begin{equation}\label{eq:new-necessary-same-ramified}
  \frac{D}{D'}v^2-\lrb{2\varpi^n\alpha_{\theta{\overline{\chi}}}\sqrt{D}-2\sqrt{\frac{D}{D'}}}v+(1-Du^2)\equiv 0\mod\varpi^{n-\lfloor \frac{l}{2}\rfloor}.
\end{equation}
This implies that for fixed $u$, the discriminant of the equation 
\begin{equation}
\Delta(u)=4\varpi^{n}\alpha_{\theta{\overline{\chi}}}\sqrt{D}\lrb{\varpi^{n}\alpha_{\theta{\overline{\chi}}}\sqrt{D}-2\sqrt{\frac{D}{D'}}}+4\frac{D}{D'}Du^2
\end{equation}
has to be a square $\mod \varpi^{n-\lfloor \frac{l}{2}\rfloor}$.
When $n-l$ is even, we can pick $u=0$ directly. Whether $\Delta(0)$ is a square is consistent with Tunnell-Saito's test.
When ${\Delta}(0)$ is indeed a square, we get two solutions of $v\mod \varpi^{ \lceil n/2\rceil}$. 
For each of these two solutions we have
\begin{equation}
\WaldsI=\frac{1}{q^{\lfloor l/2\rfloor}}.
\end{equation}

Now if $n-l$ is odd, $v_\F({\Delta}(0))=n-l$ is odd, thus ${\Delta}(0)$ can never be a square. We need to pick $u$ such that $v_\F(u)=\frac{n-l-1}{2}$ and ${\Delta}(0)+4\frac{D}{D'}Du^2$ can be of higher evaluation and a square. Whether this is possible is again consistent with Tunnell-Saito's test.
 In this case it's possible to get more or less solutions of $v\mod\varpi^{\lceil n/2\rceil}$. For each solution we have
\begin{equation}
\WaldsI=\frac{1}{q^{\lfloor l/2\rfloor}}.
\end{equation}
\subsection{Kirillov model for minimal vectors}\label{Kirillov}
Here we describe the minimal vectors explicitly in the Kirillov model. For this purpose, we choose the special intermediate subgroup $B^1=U_\BL(1)K_{\fA_2}(2n+1)$ in the case $e_\BL=1$ and $c(\pi)=4n+2$. Recall we choose $D'$ such that  
\begin{equation}\label{Eq:specialAlphaTheta}
 \alpha_\theta= \frac{1}{\varpi_\BL^{c(\theta)}}\frac{1}{\sqrt{D'}} \mapsto \frac{1}{\varpi^{c(\theta)/e_\BL}}\matrixx{0}{\frac{1}{D'}}{1}{0}.
\end{equation}
We define the intertwining operator from $\pi$ to its Whittaker model by
\begin{equation}\label{eq:3.4:IntertwiningtoWhittaker}
\varphi \mapsto W_\varphi(g)=\int\limits_{\BF}\Phi_\varphi\lrb{\matrixx{\varpi^{\lfloor c(\pi)/2\rfloor}}{0}{0}{1}\matrixx{1}{n}{0}{1}g}\psi(-n)dn.
\end{equation}
A particular minimal vector was given in the Kirillov model in \cite{HN18} under this operator.
\begin{lem}\label{lem:toricnewforminKirillov}
Up to a constant multiple, a minimal vector $\varphi_0$ is given in the Kirillov model by the following: \begin{enumerate}
\item When $c(\pi)=4n$, $\varphi_0=\Char(\varpi^{-2n}U_\BF(n)) $.
\item When $c(\pi)=2n+1$, $\varphi_0=\Char(\varpi^{-n}U_\BF(\lceil n/2 \rceil))$.
\item When $c(\pi)=4n+2$, $\varphi_0=\Char(\varpi^{-2n-1}U_\BF(n+1))$.
\end{enumerate}
\end{lem}
\begin{coro}\label{Cor:RelationNewMinimal}
The newform $\varphi_{new}$ can be related to $\varphi_0$ by the following formula
\begin{equation}
\varphi_{new}=\frac{1}{\sqrt{(q-1)q^{\lceil \frac{c(\theta)}{2 e_\BL}\rceil-1}}}\sum\limits_{x\in (\CO_\BF/\varpi^{\lceil \frac{c(\theta)}{2 e_\BL}\rceil} \CO_\BF)^\times} \pi\lrb{\matrixx{\varpi^{-c(\theta)/e_\BL}x}{0}{0}{1}}\varphi_0.
\end{equation}
Here $\varphi_0$ and $\varphi_{new}$ are both $L^2$-normalized.
\end{coro}
\begin{proof}
By the previous lemma one can uniformly write
\begin{equation}
\varphi_0=\sqrt{(q-1)q^{\lceil \frac{c(\theta)}{2 e_\BL}\rceil-1}}\Char(\varpi^{-c(\theta)/e_\BL} U_\BF(\lceil \frac{c(\theta)}{2 e_\BL}\rceil)).
\end{equation}
The coefficient comes from the $L^2$-normalization of $\varphi_0$, as $$\Vol(U_\BF(\lceil \frac{c(\theta)}{2 e_\BL}\rceil))=\frac{1}{(q-1)q^{\lceil \frac{c(\theta)}{2 e_\BL}\rceil-1}}.$$
Then one just has to use that $\varphi_{new}=\Char(\CO_\BF^\times)$ in the Kirillov model. 
\end{proof}

\section{Waldspurger's period integral using newforms}\label{newforms}
In the last section we reviewed the local Waldspurger's period integral for minimal vectors. In this section we show how to work out the local Waldspurger's period integral for newforms. Using the relation between the newform and the minimal vectors in Corollary \ref{Cor:RelationNewMinimal}, and the bilinearity of Waldspurger's period integral, we can write the integral for the newform as a sum of integrals for the minimal vectors. In the case when there is a single non-vanishing term, we get the integral very easily in Section \ref{Sec:Singleterm}, which turns out to be enough for the special example in Section \ref{Sec:Localcomputation}. We further illustrate how to evaluate the off-diagonal terms in Section \ref{Sec:testingInGeneral} in more general cases for possible future applications.


As in the introduction we denote
\begin{equation}\label{Eq:DefBetaI}
\BetaI{\varphi_1,\varphi_2}=\int\limits_{t\in \BF^\times\backslash\BE^\times}(\pi(t)\varphi_1,\varphi_2)\chi(t)dt
\end{equation}
for the embedding of $\BE$ as in \eqref{Eq:standardembedding}. Using its bilinearity and Corollary \ref{Cor:RelationNewMinimal}, we immediately have the following:

\begin{coro}\label{Cor:PeriodIntRelation}
Let $\widetilde{\varphi_{new}}=\pi\lrb{\matrixx{\varpi^{c(\theta)/e_\BL}}{0}{0}{1}}\varphi_{new}$, $\varphi_x=\pi\lrb{\matrixx{x}{0}{0}{1}}\varphi_0$ Then
\begin{equation}
\BetaI{\widetilde{\varphi_{new}},\widetilde{\varphi_{new}}}=\frac{1}{(q-1)q^{\lceil \frac{c(\theta)}{2 e_\BL}\rceil-1}}\sum\limits_{x,x'\in (\CO_\BF/\varpi^{\lceil \frac{c(\theta)}{2e_\BL}\rceil}\CO_\BF)^\times}\BetaI{\varphi_x,\varphi_{x'}}.
\end{equation}
\end{coro}

\begin{lem}\label{lem:Gl2-newform-crossterms}
\begin{enumerate}
\item Suppose that $\BetaI{\varphi_x,\varphi_x}=0$, then  $\BetaI{\varphi_{x'},\varphi_x}=\BetaI{\varphi_x,\varphi_{x'}}=0$ for any $x'$.
\item Suppose that $|\BetaI{\varphi_x,\varphi_x}|=|\BetaI{\varphi_{x'},\varphi_{x'}}|$. Then $|\BetaI{\varphi_{x},\varphi_{x'}}|= |\BetaI{\varphi_{x},\varphi_{x}}|$.
\end{enumerate}

\end{lem}
\begin{proof}
For any nontrivial functional $\mathcal{F}\in \Hom_{\BE^\times}(\pi\otimes \chi,\BC)$, we have $\BetaI{\varphi_1,\varphi_2}=C \mathcal{F}(\varphi_1)\overline{\mathcal{F}(\varphi_2)}$ for some non-zero constant $C$ independent of the test vectors, as $\dim \Hom_{\BE^\times}(\pi\otimes \chi,\BC)\leq 1$.
Then
$$|\BetaI{\varphi_x,\varphi_x}|=|C\mathcal{F}(\varphi_x)^2|,$$
$$|\BetaI{\varphi_{x'},\varphi_{x'}}|=|C\mathcal{F}(\varphi_{x'})^2|,$$
$$|\BetaI{\varphi_{x},\varphi_{x'}}|=|C\mathcal{F}(\varphi_{x})\overline{\mathcal{F}(\varphi_{x'})}|.$$
Now the results are clear.
\end{proof}

Note that the diagonal terms where $x=x'$ are already known by Section \ref{Sec:WaldsForMinimalVec}. 


\subsection{The special case}\label{Sec:Singleterm}
\begin{prop}\label{Prop:singlenonvanishingcase}
When $\BetaI{\varphi_x,\varphi_x}\neq 0$ for a single $x\in (\CO_\BF/\varpi^{\lceil \frac{c(\theta)}{2e_\BL}\rceil}\CO_\BF)^\times$, all off-diagonal terms vanish and
\begin{equation}
 \BetaI{\widetilde{\varphi_{new}},\widetilde{\varphi_{new}}}=\frac{1}{(q-1)q^{\lceil \frac{c(\theta)}{2 e_\BL}\rceil-1}}\BetaI{\varphi_x,\varphi_{x}}.
\end{equation}
\end{prop}
\begin{proof}
It follows from Lemma \ref{lem:Gl2-newform-crossterms} that all the off-diagonal terms vanish, and only a single diagonal terms is non-vanishing. The proposition now follows directly from Corollary \ref{Cor:PeriodIntRelation}.
\end{proof}
\begin{remark}
This case happens mostly when $c(\theta\chi)$ or $c(\theta\overline{\chi})\leq 1$. But there are other possibilities, as we shall see in the special example in Section \ref{Sec:Localcomputation}.
\end{remark}

\subsection{The general cases}\label{Sec:testingInGeneral}
Now we explain how to work out the off-diagonal terms in Corollary \ref{Cor:PeriodIntRelation} in more general cases. 
\begin{defn}\label{Defn:suppofint}
By the support of the local Waldspurger's period integral $\{\varphi,\varphi'\}$, we mean the set $\BE^\times\cap \Supp \Phi_{\varphi,\varphi'}$.
\end{defn}

The main idea is to get the size of the off-diagonal terms by Lemma \ref{lem:Gl2-newform-crossterms}, and to get the support of the integral in Lemma \ref{Lem:offdiagSupp}. The volume of the support of the integral is exactly the size of the integral, while the integrand is absolutely bounded by 1. This forces the integrand to be constant (with absolute value 1) on the support of the integral. Then one can easily detect this constant by looking at the value of the integrand at any point in the support of the integral.

For simplicity, however, we stay in the setting where $\BE\simeq \BL$ are ramified, $0<c(\theta{\overline{\chi}})=2l\leq 2n$. We further assume that $n-l$ is even. By Section \ref{Sec:WaldsForMinimalVec}, we can pick $u=0$, and there exists 2 solutions $v,v' \mod \varpi^{\lceil n/2\rceil}$ to \eqref{eq:new-necessary-same-ramified}, while the diagonal terms are always $\frac{1}{q^{\lfloor l/2\rfloor}}$ for these 2 solutions.

According to Lemma \ref{lem:Gl2-newform-crossterms}, we can write that $\BetaI{\varphi_{v},\varphi_{v'}}= \gamma \BetaI{\varphi_{v},\varphi_{v}}$ for some phase factor $\gamma$ with $|\gamma|=1$. Then
\begin{align}
I(\widetilde{\varphi_{new}},\chi)&=\frac{1}{(q-1)q^{\lceil \frac{c(\theta)}{2 e_\BL}\rceil-1}}(\BetaI{\varphi_{v},\varphi_{v}}+\BetaI{\varphi_{v},\varphi_{v'}}+\BetaI{\varphi_{v'},\varphi_{v}}+\BetaI{\varphi_{v'},\varphi_{v'}})\\
&=\frac{1}{(q-1)q^{\lceil \frac{c(\theta)}{2 e_\BL}\rceil-1}}\frac{1}{q^{\lfloor l/2\rfloor}}(1+\gamma)(1+\overline{\gamma}). \notag
\end{align}

To study $\gamma$,  we first study the support of the integral. Without loss of generality, we can assume that $v$, $v'$ satisfy
\begin{equation}\label{Eq:exactEqvv'}
 \frac{D}{D'}v^2-\lrb{2\varpi^n\alpha_{\theta{\overline{\chi}}}\sqrt{D}-2\sqrt{\frac{D}{D'}}}v+1=0,
\end{equation}
compared to \eqref{eq:new-necessary-same-ramified} while taking $u=0$.

Let $k=\matrixx{v}{0}{0}{1}$ and $k'=\matrixx{v'}{0}{0}{1}$.
Then for $t=\matrixx{a}{b}{bD}{a}$,
$$k'^{-1}tk=\matrixx{v v'^{-1}a }{v'^{-1}b}{vbD}{a},$$
and
\begin{equation}
\BetaI{\varphi_{v},\varphi_{v'}}=\int\limits_{\BF^\times\backslash\BE^\times} \Phi_{\varphi_0}(k'^{-1}tk)\chi(t)dt.
\end{equation}
\begin{lem}\label{Lem:offdiagSupp}
For the setting as above, we have $vv'D=D'$, and $v_\F(\frac{v}{v'}-1)= \frac{n-l}{2}$.
In particular the support of the integral 
is $v_\F(b)=0, v_\F(a)\geq  \lceil\frac{l+1}{2}\rceil$.
\end{lem}
\begin{proof}
According to \eqref{Eq:exactEqvv'}, $v$ and $v'$ satisfy 
$$v v'=\frac{D'}{D},\text{ \ \ } v+v'=2(\varpi^{n}\alpha_{\theta{\overline{\chi}}}\sqrt{D'}-1)\sqrt{\frac{D'}{D}}$$
so the first result is direct.
For the second result, we note that
\begin{equation}
\lrb{\frac{v}{v'}-1}^2=\frac{(v+v')^2-4vv'}{v'^2}=\frac{D'}{D}\frac{4\varpi^{n}\alpha_{\theta{\overline{\chi}}}\sqrt{D'}(\varpi^{n}\alpha_{\theta{\overline{\chi}}}\sqrt{D'}-2)}{v'^2}.
\end{equation}
Thus $v_\F((\frac{v}{v'}-1)^2)=n-l$ and $v_\F(\frac{v}{v'}-1)= \frac{n-l}{2}$. 
Now for $k'^{-1}tk\in J=\BL^\times K_{\fA_{e_\BL}}(n)$, there are two parts to consider: either  $v_\F(b)=0$, $v_\F(a)>0$ or $v_\F(a)=0$, $v_\F(b)\geq 0$. 
In the first case, since $$\varpi_\BL \CB_{e_\BL} ^n=\matrixx{\varpi^{ \lceil \frac{n+1}{2}\rceil}\CO_\BF}{\varpi^{ \lfloor \frac{n+1}{2}\rfloor}\CO_\BF}{\varpi^{ \lfloor \frac{n+1}{2}\rfloor+1}\CO_\BF}{\varpi^{ \lceil \frac{n+1}{2}\rceil}\CO_\BF}$$ we must have
\begin{equation}
a(\frac{v}{v'}-1)\equiv 0 \mod \varpi^{\lceil \frac{n+1}{2}\rceil}, \text{\ \ }vbD-v'^{-1}bD'\equiv 0 \mod \varpi^{\lfloor\frac{n+1}{2}\rfloor +1}.
\end{equation}
The second equation is automatic as $v v'=\frac{D'}{D}$. From the first equation and the computation for $v_\F(\frac{v}{v'}-1)$ above, we get $v_\F(a)\geq \lceil\frac{l+1}{2}\rceil$.
Using similar argument, one can easily show that it is impossible for $k'^{-1}tk\in \BL^\times K_{\fA_{e_\BL}}(n)$ when $v_\F(a)=0$, $v_\F(b)\geq 0$.
\end{proof}
\begin{prop}\label{Prop:TestingonNewform}
Suppose that $\BE\simeq\BL$ are ramified, $c(\theta)=2n$, $w_\pi=1$, $\epsilon(\pi_\BE\times\chi)=1$, and $0<c(\theta\overline{\chi})=2l\leq 2n$ with $n-l$ even. Then
$$I(\widetilde{\varphi_{new}},\chi)=\frac{1}{(q-1)q^{\lceil \frac{c(\theta)}{2 e_\BL}\rceil-1}}\frac{1}{q^{\lfloor l/2\rfloor}} (1+\theta\chi(\sqrt{D}))^2    .$$
In particular it is either $0$ or asymptotically $\frac{1}{C(\pi\times\pi_{\chi})^{1/4}}$.
\end{prop}
\begin{proof}
We already know that $|\BetaI{\varphi_{v},\varphi_{v'}}|=\frac{1}{q^{\lfloor l/2\rfloor}} $. By the previous lemma the support of the integral also has volume $\frac{1}{q^{\lfloor l/2\rfloor}}$, while the integrand  $|\pair{\pi(t)\varphi_v,\varphi_{v'}}\chi(t)|\leq 1$. Then  $\pair{\pi(t)\varphi_v,\varphi_{v'}}\chi(t)$ must be some constant $\gamma$ on the whole support with $|\gamma|=1$. To detect this constant we just have to take $t=\matrixx{0}{1}{D}{0}$. Then
\begin{equation}
\gamma=\Phi_{\varphi_0}\lrb{\matrixx{0}{v'^{-1}}{vD}{0}}\chi(\sqrt{D})=\theta(v'^{-1}\sqrt{D'})\chi(\sqrt{D})=\theta\chi(\sqrt{D}).
\end{equation}
In the last equality we have used $\theta|_{\BF^\times}=1$ so that $$\theta(v'^{-1}\sqrt{D'})=\theta\lrb{v'^{-1}\frac{\sqrt{D'}}{\sqrt{D}}\sqrt{D}}=\theta(\sqrt{D}).$$ Note that $\theta\chi(\sqrt{D})=\pm 1$.
The last statement is easy to check.
\end{proof}

\section{A special example from arithmetic geometry}\label{Sec:Localcomputation}
Now we specialize to the case required in the proof of \cite[Theorem 4.3]{HSY}. 
We shall first review the global arithmetic setting, and use the subscripts to indicate the local components.

For an arbitrary nonzero integer $n$, let $E_n$ be the elliptic curve defined by the affine equation $x^3+y^3=n$. Then $E_n$ has complex multiplication by the field $K=\BQ(\sqrt{-3})$. For a prime $p\equiv 4,7\mod 9$, the elliptic curve $E_p$ is closely related to the well-known Sylvester conjecture. Let $\pi$ be the automorphic representation of $\GL_2(\BQ)$ corresponding to $E_9$ and $\pi_3$ the $3$-adic part of $\pi$. Then $c(\pi_3)=5$. 
Assume that $f_3$ is the standard newform of $\pi_3$. Let $\chi:\Gal(\bar{K}/K)\ra \CO_K^\times$ be the character given by $\chi(\sigma)=(\sqrt[3]{3p})^{\sigma-1}$. We also view $\chi$ as a Hecke character on $\BA_K$ by the Artin map.
Then $c(\chi_3)=4$. 
The following normalized Waldspurger's period integral
\begin{equation}\label{wpi}
\beta^0_3(f_3, f_3)
=\int\limits_{t\in \BQ_3^\times\backslash K_3^\times}\frac{ (\pi(t)f_3,f_3)}{(f_3,f_3)}\chi_3(t)dt
\end{equation}
appears in the proof of the explicit Gross-Zagier formula for $E_p$ in \cite{HSY}. Here $K$ is embedded into $\rm{M}_2(\BQ)$ as in \cite[Section 2]{HSY} which linearly extends the following map:
\begin{equation}\label{matrix}
\sqrt{-3}\mapsto \matrixx{4p+17+72/p}{-8p/9-4-18/p}{18p+72+288/p}{-4p-17-72/p}=:\matrixx{a}{3^{-2}b}{3^3c}{-a}
\end{equation}
with $3||a$ if $p\equiv 4\mod 9$, $9||a$ if $p\equiv 7\mod 9$, $b\equiv p\mod 9$ and $c\equiv -1\mod 9$.  Then $\Nm(\sqrt{-3})=-a^2-3bc=3$.
 Note here that $K$ is embedded differently from the fixed embedding we have been using in \eqref{Eq:standardembedding}. We choose the notation $\beta_3^0$  to reflect this difference, and to be consistent with the notation in \cite{HSY}. We shall  work out the relation between these two embeddings later on.


To apply the results in the previous section to compute (\ref{wpi}), we take $\varpi=3=q$, $D=-3$, $K_3\simeq\BE\simeq\BL\simeq\BQ(\sqrt{-3})_3$, $c(\theta_3)=c(\chi_3)=4$, $n=2$.
By Lemma \ref{lem:toricnewforminKirillov} we have the minimal vector $\varphi_0=\Char(\varpi^{-2}U_\BF(1))$ in the Kirillov model. As seen from previous sections, more accurate information is needed. 
\subsection{Local characters associated to the arithmetic information}
First of all we make use of the arithmetic information to give the local characters explicitly. Recall that $K=\BQ(\sqrt{-3})$ is an imaginary quadratic field and $\CO_K=\BZ[\omega]$ is its ring of integers with $\omega=\frac{-1+\sqrt{-3}}{2}$. Let $\Theta:K^\times\backslash \BA_K^\times\ra \BC^\times$ be the unitary Hecke character associated to the base-changed CM elliptic curve ${E_9}_{/K}$. 
Then $\Theta$ has conductor $9\CO_K$. For any place $v$ of $K$, let $\Theta_v$ be the local component of $\Theta$ at the place $v$. Then $\Theta_v$ is the character used to construct the local Weil-Deligne representation in Theorem \ref{LLC}. We denote $\Theta_3$ the 3-part of $\Theta$. Then $\pi_3$ is the local representation of $\GL_2(\BQ_3)$ corresponding to $\Theta_3$. Note
\[\CO_{K,3}^\times/(1+9\CO_{K,3})\simeq \langle \pm 1\rangle ^{\BZ/2\BZ} \times\langle 1+\sqrt{-3}\rangle^{\BZ/3\BZ}\times\langle 1-\sqrt{-3}\rangle^{\BZ/3\BZ}\times \langle 1+3\sqrt{-3}\rangle^{\BZ/3\BZ}.\]
\begin{lem}\label{thetavalue}
We have $c(\Theta_3)=4$. Its values are
given explicitly by
\[\Theta_3(-1)=-1,\quad \Theta_3(1+\sqrt{-3})=\frac{-1-\sqrt{-3}}{2},\ \ \Theta_3(\sqrt{-3})=i,\]
\[\Theta_3(1-\sqrt{-3})=\frac{-1+\sqrt{-3}}{2},\quad \Theta_3(1+3\sqrt{-3})=\frac{-1+\sqrt{-3}}{2}.\]
\[\]
\end{lem}
\begin{proof}
It is well-known that $\Theta_\infty(x)=\frac{||x||}{x}$, (see for example \cite[Chapter II, Theorem 9.2]{Silvermanbook2} and note that we normalize it to make it unitary) and $\Theta$ is unramfied when $3\nmid v$. Note
\[\Theta_\infty(-1)\Theta_3(-1)=1,\quad \Theta_\infty(-1)=-1.\]
So $\Theta_3(-1)=-1$.

Let $\fp=(a)$ be a prime of $K$  relatively prime to $6$, with the unique generator $a\equiv 2\mod 3$.
By \cite[Chapter II, Example 10.6]{Silvermanbook2}, we have 
\[\Theta(\fp)=-\RN_\fp^{-1/2}\ov{\left(\frac{-3}{a}\right)}_6a.\]
Where $\left(\frac{\cdot}{a}\right)_6$ is the sixth power residue symbol and $\RN_\fp$ is the norm of $\fp$.
If $\fp=(5)$, then
\[\Theta_5(5)=-\ov{\left(\frac{-3}{5}\right)}_6=-1.\]
By
\[\Theta_\infty(10)\Theta_2(10)\Theta_3(10)\Theta_5(10)=1,\]
we get $\Theta_2(2)=-1.$ Since $\Theta_2$ is unramified and $1+\sqrt{-3}$ is another uniformizer of $(2)$, we see $\Theta_2(1\pm\sqrt{-3})=\Theta_2(2)=-1$.
By
\[\Theta_\infty(1\pm\sqrt{-3})\Theta_2(1\pm\sqrt{-3})\Theta_3(1\pm\sqrt{-3})=1,\quad \Theta_\infty(1\pm\sqrt{-3})=\frac{2}{1\pm\sqrt{-3}},\]
we get
\[\Theta_3(1\pm\sqrt{-3})=\frac{-1\mp\sqrt{-3}}{2}.\]
The other evaluations can be obtained in a similar way.


\end{proof}

The local  character  $\chi_3$ has conductor $\BZ_3^\times(1+9\CO_{K,3})$, and hence it is in fact a character of the quotient group $\CO_{K,3}^\times/\BZ_3^\times(1+9\CO_{K,3})$. Note that
\[\CO_{K,3}^\times/\BZ_3^\times(1+9\CO_{K,3})\simeq \langle 1+\sqrt{-3}\rangle^{\BZ/3\BZ}\times\langle 1+3\sqrt{-3}\rangle^{\BZ/3\BZ}.\]
We have the following, 
\begin{lem} \label{chi}
$c(\chi_3)=4$ and $\chi_3|_{\BQ_3^\times}=1$. Its values are given explicitly by the following table:
\begin{center}
\begin{tabular}{|c|c|c|c|c|c|}
\hline
$p\mod 9$&$\chi_3(1+\sqrt{-3})$&$\chi_3(1+3\sqrt{-3})$&$\chi_3(\sqrt{-3})$\\
\hline
$4$&$\omega$&$\omega$&$1$\\
\hline
$7$&$\omega^2$&$\omega$&$1$\\
\hline
\end{tabular}
\end{center}
\end{lem}

\begin{proof}
This follows directly from the explicit local class field theory.
We note that all the elements in $1+9\CO_{K,3}$ is a cube in $K_3$. 
Hence for any $t\in K_3^\times$,
\[\chi_3(t)=\left(\sqrt[3]{3p}\right)^{\sigma_t-1}=\lrb{\frac{t,3p}{K_3;3}}=\left\{\begin{aligned}\lrb{\frac{t,12}{K_3;3}},&\quad p\equiv 4\mod 9;\\ \lrb{\frac{t,21}{K_3;3}},&\quad p\equiv 7\mod 9.\\\end{aligned}\right.\]
Recall $\sigma_t$ is the  image of $t$ under the the Artin map, and $\lrb{\frac{\cdot,\cdot}{K_3;3}}$ denotes the $3$-rd Hilbert symbol over $K_3^\times$ as in \cite{HSY}. Using the local and global principal, it is straight-forward to compute the values of $\chi_3$ as in the above table. For example,  
\[\chi_3(1+\sqrt{-3})=\lrb{\frac{1+\sqrt{-3},3}{K_3;3}}\lrb{\frac{1+\sqrt{-3},4}{K_3;3}}=\lrb{\frac{1+\sqrt{-3},3}{K_2;3}}^{-1}\lrb{\frac{1+\sqrt{-3},4}{K_2;3}}^{-1}=\omega,\]
where the second equality uses the fact that the product of the Hilbert symbol at all places equals $1$ and the Hilbert symbol we considered is unramified outside the place $2$ and $3$,  the third
equality uses the formula in \cite[Chapter V, Proposition 3.4]{Neukirchbook1}. 
\end{proof}

\subsection{Local period integral}

From Section \ref{Sec:WaldsForMinimalVec}, we see that the test vector issue for Waldspurger's period integral is closely related to $c(\theta_3\ov\chi_3)$ or $c(\theta_3\chi_3)$ and some further details like $\alpha_{\theta_3\ov\chi_3}$. We can work out these  by using Lemma \ref{thetavalue}, \ref{chi}, and the relation between $\theta_3$ and $\Theta_3$ as in Theorem \ref{LLC}.
\begin{coro}
\label{thetachivalue}
If $p\equiv 4\mod 9$,
the local character $\Theta_3\ov\chi_3$ is given explicitly by
\[\Theta_3\ov\chi_3(-1)=-1,\quad \Theta_3\ov\chi_3(1+\sqrt{-3})=\omega,\]
\[\Theta_3\ov\chi_3(1-\sqrt{-3})=\omega^2,\quad \Theta_3\ov\chi_3(1+3\sqrt{-3})=1,\quad \Theta_3\ov\chi_3(\sqrt{-3})=i.\]

If $p\equiv 7\mod 9$,
the local character $\Theta_3\ov\chi_3$ is given explicitly by
\[\Theta_3\ov\chi_3(-1)=-1,\quad \Theta_3\ov\chi_3(1+\sqrt{-3})=1,\]
\[\Theta_3\ov\chi_3(1-\sqrt{-3})=1,\quad \Theta_3\ov\chi_3(1+3\sqrt{-3})=1,\quad \Theta_3\ov\chi_3(\sqrt{-3})=i.\]
\end{coro}
Now we can prove the following key Lemma in our special case.
\begin{lem}\label{Cor:AllnecessaryFormulationForThetaChi}
When $p\equiv 7\mod 9$, we have $\theta_3\ov\chi_3=1$. When $p\equiv 4\mod 9$, we have $c(\theta_3\ov\chi_3)=2$ and $\alpha_{\theta_3\ov\chi_3}=\frac{1}{3\sqrt{-3}}$.
\end{lem}
\begin{proof}
Let $\psi_3$ be the additive character such that $\psi_3(x)=e^{2\pi i \iota(x)}$ where $\iota:\BQ_3\rightarrow \BQ_3/\BZ_3 \subset \BQ/\BZ$ is the map given by $x\mapsto -x\mod \BZ_3$ which is compatible with the choice in \cite{CST14}. Let $\psi_{K_3}(x)=\psi_3\circ \Tr_{K_3/\BQ_3}(x)$, be the additive character of $K_3$. 

Recall that $\alpha_{\Theta_3}$ is the number associated to $\Theta_3$ as in Lemma \ref{Lem:DualLiealgForChar} so that
\[\Theta_3(1+x)=\psi_{K_3}(\alpha_{\Theta_3} x), \]
for any $x$ satisfying $v_{K_3}(x)\geq c(\Theta_3)/2=2$. By the definition of $\psi_{K_3}$ and Lemma \ref{thetavalue}, we know that $\alpha_{\Theta_3}=\frac{1}{9\sqrt{-3}}$.  Now let $\eta_3$ be the quadratic character associated to the quadratic field extension $K_3/\BQ_3$.
Then by \cite[Proposition 34.3]{BushnellHenniart:06a}, $\lambda_{K_3/\BQ_3}(\psi')=\tau(\eta_3,\psi_3')/\sqrt{3}=-i$, here $\tau(\eta_3,\psi_3')$ is the Gauss sum and $\psi_3'(x)=\psi_3(\frac{x}{3})$ is the additive character of level one. By \cite[Lemma 5.1]{Langlands}, $\lambda_{K_3/\BQ_3}(\psi_3)=\eta_3(3)\lambda_{K_3/\BQ_3}(\psi_3')=-i$. Then $\Delta_{\theta_3}$ is the unique level one character of $K_3$ such that $\Delta_{\theta_3}|_{\BZ_3^\times}=\eta_3$ and
\[\Delta_{\theta_3}(\sqrt{-3})=\eta((\sqrt{-3})^3\alpha_{\Theta_3})\lambda_{K_3/\BQ_3}(\psi_3)^{3}=-i.\]
Recall that $\theta_3=\Theta_3\Delta_{\theta_3}$. Then by Corollary \ref{thetachivalue} we can easily check that: 

\begin{enumerate}
\item If $p\equiv 7\mod 9$, $\theta_3\ov\chi_3$ is the trivial character. 
\item If $p\equiv 4\mod 9$, $\theta_3\ov\chi_3$ is of level 2 and by definition we can choose $\alpha_{\theta_3\ov\chi_3}=\frac{1}{3\sqrt{-3}}$.
\end{enumerate}
\end{proof}

\begin{prop}\label{Prop:TestingForNew}
Suppose $\Vol(\BZ_3^\times\backslash\CO_{K,3}^\times)=1$ so that $\Vol(\BQ_3^\times\backslash K_3^\times)=2$.
For $f_3$ being the newform, $K$ being embedded as in \eqref{matrix} and $\theta_3,\chi_3$ as given above, we have
\begin{equation}
 \beta^0_3(f_3,f_3)=\begin{cases}
 1, &\text{\ if }p\equiv 7\mod 9\\
 1/2, &\text{\ if }p\equiv 4\mod 9.
 \end{cases}
\end{equation}
\end{prop}

\begin{proof}
We may assume $f_3$ to be $L^2$-normalized. To evaluate $f_3$ for the embedding in \eqref{matrix} is equivalent to use the standard embedding \eqref{Eq:standardembedding} of $\BE$ and use different translate of the newform.
In particular the embedding in \eqref{matrix} is conjugate to the standard embedding by the following. 
\begin{equation}
    \matrixx{a}{3^{-2}b}{3^3c}{-a}=\matrixx{-9c}{a/3}{0}{1}^{-1}\matrixx{0}{1}{D}{0}\matrixx{-9c}{a/3}{0}{1}.
\end{equation}
Thus  we have
\begin{align}
\beta^0_3(f_3,f_3)&=\int\limits_{\BF^\times\backslash \BE^\times}\lrb{\pi_3\lrb{\matrixx{-9c}{a/3}{0}{1}^{-1}t\matrixx{-9c}{a/3}{0}{1}}f_3,f_3}\chi(t)dt\\
&=\int\limits_{\BF^\times\backslash \BE^\times}\lrb{\pi_3\lrb{t\matrixx{-9c}{a/3}{0}{1}}f_3,\pi_3\lrb{\matrixx{-9c}{a/3}{0}{1}}f_3}\chi(t)dt\notag,
\end{align}
which is by definition $\BetaI{\pi_3\lrb{\matrixx{-9c}{a/3}{0}{1}}f_3,\pi_3\lrb{\matrixx{-9c}{a/3}{0}{1}}f_3}$ for the bilinear pairing as in \eqref{Eq:DefBetaI} and the standard embedding as in \eqref{Eq:standardembedding}.
Note that by Corollary \ref{Cor:RelationNewMinimal}
$$\pi_3\lrb{\matrixx{-9c}{a/3}{0}{1}}f_3
=\frac{1}{\sqrt{2}}\sum\limits_{x\in (O_v/\varpi O_v)^\times}\pi_3\lrb{\matrixx{1}{a/3}{0}{1}\matrixx{x}{0}{0}{1}}\varphi_0$$
where $\varphi_0$ is the minimal test vector.

Now there are two cases. When $p\equiv 7\mod 9$, $9|| a$ and the action of $\matrixx{1}{a/3}{0}{1}$ on $\varphi_x=\pi_3\lrb{\matrixx{x}{0}{0}{1}}\varphi_0$ is by a simple character. By the $l=0$ case in Section \ref{Sec:WaldsForMinimalVec}, we have a unique $x\mod\varpi$ for which $\BetaI{\varphi_x,\varphi_x}$ is nonvanishing. According to Proposition \ref{Prop:singlenonvanishingcase}
 there are no off-diagonal terms, and we have
\begin{equation}
\beta^0_3\lrb{f_3,f_3 }=\frac{1}{(q-1)q^{\lceil \frac{c(\theta_3)}{2 e_\BL}\rceil-1}} \BetaI{\varphi_x,\varphi_x}=\frac{1}{2}\cdot 2=1.
\end{equation}

When $p\equiv 4\mod 9$, we have $3||a$ and $u=a/3$. By Lemma \ref{Cor:AllnecessaryFormulationForThetaChi}, this is the case $l=1$ and $n-l=1$ is odd. By the choice in Section \ref{Sec:WaldsForMinimalVec}, 
$$D'=\frac{1}{\alpha_{\theta_3}^2\varpi_\BL^{2c(\theta_3)}}=-3.$$

By Lemma \ref{Cor:AllnecessaryFormulationForThetaChi}, $\alpha_{\theta_3\chi_3^{-1}}=\frac{1}{3\sqrt{-3}}$ in this case, and we have 
\begin{align}
\Delta(u)&= 4\varpi^{n}\alpha_{\theta_3\ov\chi_3}\sqrt{D}\lrb{\varpi^{n}\alpha_{\theta_3\ov\chi_3}\sqrt{D}-2\sqrt{\frac{D}{D'}}}+4\frac{D}{D'}Du^2 \\
&\equiv 4\cdot 9\cdot \frac{1}{3\sqrt{-3}} \cdot \sqrt{-3}\cdot (-2)+4 \cdot (-3)\frac{a^2}{9}\mod{\varpi^2}\notag\\
&\equiv  -8\cdot 3-4\cdot 3 \mod{\varpi^2}\notag\\
&\equiv 0 \mod{\varpi^2}\notag.
\end{align}
$\Delta(u)$ is indeed  congruent to a square. Then we can get a unique solution of $v\mod \varpi$ from \eqref{eq:new-necessary-same-ramified}, and again by Proposition \ref{Prop:singlenonvanishingcase},
\begin{equation}
 \beta^0_3\lrb{f_3,f_3}=\frac{1}{(q-1)q^{\lceil \frac{c(\theta_3)}{2 e_\BL}\rceil-1}}\frac{1}{q^{\lfloor l/2\rfloor}}=\frac{1}{2}.
\end{equation}
\end{proof}

Let $f'$ be the admissible test vector of $(\pi, \chi)$ which is as defined in \cite[Definition 1.4]{CST14}. By definition, the $3$-adic part $f_3'$ is $\chi_3^{-1}$-eigen under the action of $K_3^\times$. The following corollary is most directly used in \cite{HSY}.
\begin{coro}\label{ration}
For the admissible test vector $f'_3$ and the newform $f_3$ we have 
\[\frac{\beta_3^0(f'_3,f'_3)}{\beta_3^0(f_3,f_3)}=\begin{cases}
 2, &\text{\ if }p\equiv 7\mod 9,\\
 4, &\text{\ if }p\equiv 4\mod 9.
 \end{cases}\]
\end{coro}
\begin{proof}
Keep the normalization of the volumes in Proposition \ref{Prop:TestingForNew}. By definition of $f'$, we have $\beta_3^0(f'_3,f'_3)=\Vol(\BQ_3^\times\backslash K_3^\times)=2$. Then the corollary follows from Proposition \ref{Prop:TestingForNew}.
\end{proof}

\bibliographystyle{alpha}
\bibliography{reference}
\end{document}